\documentclass{amsart}
\usepackage{amsthm,amsfonts,amsmath,amssymb,latexsym,epsfig,eucal,enumitem}

\newtheorem{theorem}{Theorem}

\newtheorem{corollary}[theorem]{Corollary}

\def\mb#1{{\mathbb #1}}

\def\<{\langle}                                    
\def\>{\rangle}

\begin{document}

\title[Riemannian metric representatives of the Stiefel-Whitney classes] 
{Riemannian metric representatives of the Stiefel-Whitney classes}
\author{Santiago R. Simanca}
\thanks{Supported by the Simons Foundation Visiting Professorship award 
number 657746.}
\address{Department of Mathematics, Courant Institute of Mathematical Sciences,
251 Mercer St., New York, NY 10012} 
\email{srs2@cims.nyu.edu. {\it Current:} srsimanca@gmail.com.}

\begin{abstract} 
If $M$ is a closed manifold, and $K$ is a smooth triangulation of $M$,  
Whitney proved that all of the Stiefel-Whitney classes are specified as 
cochains  on the dual cell complex $(K')^*$ assigning the value 
$1$ mod $2$ to each dual cell. We provide the pair $(M,K)$ with an arbitrary 
Riemannian metric $g$, and use Whitney's criteria to show that there are 
associated representatives
of all the Stiefel-Whitney classes $w_1(M), \ldots , w_n(M)$. 
The representative of $w_1(M)$ is determined by $\det{g_{ij}}$, the $g_{ij}$s 
computed in a frame that is locally defined at each dual $1$-cell; the 
representatives of the even classes $w_{2k}(M)$ are determined by the 
Chern-Gauss-Bonnet density $2k$-form of locally defined totally 
geodesic oriented $2k$ manifolds with boundary associated to each 
dual $2k$-cell; and
the representatives of the odd classes $w_{2k+1}(M)$ are determined by the 
hypersurface area form of the boundary sphere of a locally defined 
totally geodesic oriented $(2k+1)$ manifold with boundary associated to each
dual $(2k+1)$-cell. If $(M,J,g)$ is Hermitian, we prove that
the metric representative of $w_{2k}(M)$ so obtained is the $\mb{Z}/2$
reduction of the $k$th Chern class $c_k(M,J)$ induced by the coefficient
homomorphism, and that the metric
representative of any odd degree class $w_{2k+1}(M)$ so obtained is 
trivial in cohomology. 
\end{abstract}


\maketitle 

\section{The Stiefel-Whitney classes}
A smooth triangulation $K$ of a closed manifold $M$ suffices to specify
all of the Stiefel-Whitney classes of $M$. For if $K'$ is the barycentric
subdivision of $K$, and $(K')^*$ is the cell complex dual to $K'$, then 
$w_i(M)$ is represented by the mod $2$ cochain that assigns $1$ to each
dual $i$-cell. This remarkably beautiful (and simple) characterization of these 
cohomology classes was announced by Whitney in \cite{whit}, but 
his proof never appeared in print. In 1970, Cheeger provided
the only known written proof of such to date \cite{chee}. 
In 1971, Halperin and Toledo gave a proof by characterizing the 
Stiefel-Whitney homology classes \cite{hato}, an equivalent 
version that had been conjectured much earlier by 
Stiefel \cite{stie}, and which was the starting point of Whitney, who 
dualized the statement to the cohomology classes.
Sullivan used the homology version to define Stiefel-Whitney 
classes for more general spaces \cite{sull}.

%{\it Added note}: 
(Upon the posting of our article in the \verb+arXiv+,
Professor Harvey brought up to our attention his paper \cite{hazw}, coauthored 
with J. Zweck. In this publication, it is proved that certain currents, which 
are associated to sections of a real vector bundle satisfying a mild 
condition, represent the Stiefel-Whitney classes of the bundle.)

When we endow a manifold with a Riemannian metric, the manifold acquires a
geometry that we attempt to understand for various reasons, one of which
is to be able to try to read off from it the topology.
If a canonical metric exists, and determining when this is the case 
is a very worthwhile, and usually hard problem in its own right, we choose 
this metric to make the geometry of the manifold the best, so 
that we can then read off its topology with ease. But the topology of
the manifold remains fixed, no matter its geometric shape, and any choice
of the metric must lead to the same readings of it, if at all.

We provide the triangulated manifold $(M,K)$ with an arbitrary 
Riemannian metric $g$. A dual $i$-cell in $(K')^*$ corresponds to a 
unique $(n-i)$-simplex $\sigma_{n-i}$ in $K$ that is a face of an $n$-simplex
$\sigma_n$ in $K$, and a sequence of simplices $\sigma_j$, $j=n-i, \ldots, n$,
where $\sigma_{j+1}$ succeeds $\sigma_j$, $j=n-i, \ldots, n-1$, in the
natural order of the simplices of $K'$.     
We firstly use the metric to fix compatibly a local positive orientation for 
all of the simplices in this chain, and relying on Whitney's characterization 
of the Stiefel-Whitney classes, we then prove the following: 
\begin{itemize}
\item Given a $1$-cell corresponding to a simplex $\sigma_{n-1}$,
we find an orthonormal frame of $T\sigma_n\mid_{\stackrel{\circ}{\sigma}_{n-1}}=
T{\stackrel{\circ}{\sigma}_{n-1}}\oplus \nu({\stackrel{\circ}{\sigma}_{n-1}})$
that fixes compatibly the local positive orientation of the simplices, and 
that we use to compute the components of the metric. Then 
the $1$-cochain defined by evaluating the function $\det{g_{ij}}$ at the 
barycenter of the simplex represents the first Stiefel-Whitney class $w_1(M)$.
\item Given a $2$-cell corresponding to a simplex $\sigma_{n-2}$, the fiber
of the normal bundle $\nu(\stackrel{\circ}{\sigma}_{n-2})$ of
$\stackrel{\circ}{\sigma}_{n-2}\subset M$ through the barycenter of the
simplex determines a local, totally geodesic, positively oriented, smooth 
$2$-disk with boundary. 
If $r$ is the intrinsic Ricci tensor of this 
disk, the limit as $\varepsilon \searrow 0$ of the evaluation of 
$\frac{1}{2\pi\varepsilon^2} r$ as a 2-form over 
the closure of the $2$-disk defines a $2$-cochain that represents the 
second Stiefel-Whitney class $w_2(M)$.
\item Given a $3$-cell corresponding to a simplex $\sigma_{n-3}$, 
the fiber of the normal bundle $\nu(\stackrel{\circ}{\sigma}_{n-3})$ of
$\stackrel{\circ}{\sigma}_{n-3}\subset M$ through the barycenter of the
simplex determines a local, totally geodesic, positively oriented smooth 
$3$-disk with boundary a $2$-sphere. If $d\sigma$ is the area form on this 
sphere, and $\varepsilon$ is
its radius, the limit as $\varepsilon \searrow 0$ of the evaluation of  
$d\sigma/4\pi \varepsilon^2$ over the $2$-sphere cycle
defines a $3$-cochain that represents the third
Stiefel-Whitney class $w_3(M)$.
\end{itemize}

It becomes clear how to extend the last argument to produce representatives
of the odd degree classes $w_{2k+1}(M)$, $k>1$. The extension of the 
second case to finding representatives of all the even degree classes 
$w_{2k}(M)$, $k>1$, is a bit more elaborate, but becomes clear also once we
notice that the role that the intrinsic Ricci tensor $r$ plays in the 
argument is exactly that played by the  
curvature $2$-form of the Levi-Civita connection of the intrinsic metric
on the totally geodesic oriented disk. The extension is then accomplished by 
using the $k$-fold product of this local curvature form in the role that $r$
plays in the case of $w_2$. Up to a suitable constant, 
this is nothing but the locally defined
density of the Chern-Gauss-Bonnet theorem for $2k$-manifolds with boundary.

When $M$ carries a complex structure $J$, and $(M,J,g)$ is Hermitian, by the
definition of the metric representative of $w_{2k}(M)$,  we see that this
agrees with the mod $2$ reduction of the $k$th Chern class $c_k(M,J)$ 
induced by the homomorphism $\mb{Z} \rightarrow \mb{Z}/2$. By studying the 
relationship between the local positive orientations of the dual cells, and 
the orientation of $M$ determined by $J$, we see also that the metric 
representative of $w_{2k+1}(M)$ is trivial in cohomology.   

\section{Smooth triangulations, barycentric subdivisions,  and dual cell 
complexes: The class $w_0(M)$}\label{smt}
We let $M=M^n$ be a smooth closed $n$-dimensional manifold, and $K$     
be a locally finite, smooth, simplicial complex triangulation of $M$, so the 
underlying polytope of $K$ is $M$. We recall quickly the notions of
barycentric subdivision $K'$ of $K$, and dual cell complex
$(K')^*$. We refer the reader to \cite{mubo} for details.
 
The barycentric subdivision $K'$ of any triangulation $K$ is a naturally 
oriented simplicial complex. The simplices of $K'$ are all of the form 
$$
\hat{\sigma}_{i_1}\hat{\sigma}_{i_2} \ldots \hat{\sigma}_{i_k}
$$
where $\sigma_{i_1}, \ldots, \sigma_{i_k}$ are simplices of $K$ such that 
$\sigma_{i_1}\succ \sigma_{i_2} \succ \cdots \succ \sigma_{i_k}$.  
Here, for any simplex $\sigma$ of $K$, we denote by 
$\hat{\sigma}$ its barycenter in $K'$, and $\sigma_i \succ \sigma_j$ means
 that the simplex $\sigma_j$ is a proper face of the simplex $\sigma_i$. 
The vertices of $K'$ are ordered by decreasing dimension of the 
simplices of the triangulation $K$ of which they are the barycenters. This 
ordering induces a linear ordering of the vertices of each simplex of $K'$. 

Given a simplex $\sigma $ in $K$, the union of all open simplices of $K'$ of 
which $\hat{\sigma}$ is the initial vertex is the interior 
$\stackrel{\circ}{\sigma}$ of $\sigma$. The block $D(\sigma)$ dual to 
$\sigma$ is the union 
of all the open simplices of $K'$ of which $\hat{\sigma}$ is the final vertex. 
The closed block $\overline{D}(\sigma)$ is the closure of $D(\sigma)$,  
and coincides with the union of all simplices of $K'$
of which $\hat{\sigma}$ is the final vertex. It is the polytope of a 
subcomplex of $K'$. 
We let $\dot{D}(\sigma)=\overline{D}(\sigma) \setminus D(\sigma)$. 

The collection $\{ D(\sigma)\}_{\sigma \in K}$ of all dual blocks is 
pairwise disjoint, and their union equals $M$. 
This collection is therefore a cell complex that we denote by $(K')^*$.
An $i$-cell in $(K')^*$ is determined by an $(n-i)$-simplex $\sigma_{n-i}$ 
in $K$, and is the interior of a simplex of the form $e^i_{\sigma_{n-i}}=[
\hat{\sigma}_n, \ldots, \hat{\sigma}_{n-i}]\subset \overline{D}(\sigma_{n-i})$
in $K'$, where $\sigma_n \succ \cdots \succ \sigma_{n-i}$.  
Notice that if $\sigma $ is a $k$-simplex of $K$, then
\cite[Theorem 64.1]{mubo}:
\begin{enumerate}
\item $\overline{D}(\sigma)$ is the polytope of a subcomplex of $K'$ of 
dimension $n-k$. 
\item $\dot{D}(\sigma)$ is the union of the blocks $D(\tau)$ for which $\tau$
has $\sigma$ as a proper face. These blocks have dimension less than 
$n-k$.
\item \label{th} If $H_i(K,K \setminus \hat{\sigma})\cong \mb{Z}$ for $i=n$ and
vanishes otherwise, then $(\overline{D}(\sigma),\dot{D}(\sigma))$ has the 
homology of an $(n-k)$-cell modulo its boundary.
\end{enumerate}

By \cite[Theorems 1.4 and 1.5]{mu}, any triangulation $K$ is isotopic to a 
smooth triangulation whose dual cells form a smooth cell decomposition of
$M$. Thus, our assumption that $K$ has this property does not restrict 
the generality of our work, and so in (\ref{th}) above, we actually have that 
$(\overline{D}(\sigma),\dot{D}(\sigma))$  is topologically an $(n-k)$-cell 
modulo its boundary.

For the convenience of the exposition, we specify the following cases:
\begin{enumerate}[label=(\alph*)] 
\item \label{ce0} If $d_0$ is a dual $0$-block, there exists an  
$n$-simplex $\sigma_n$ in $K$ such that 
$$
\overline{d}_0 = \overline{D}(\sigma_n)=\hat{\sigma}_n \, .
$$
Thus, any $0$-cell $e^0_{\sigma_{n}}$ in $(K')^*$ is determined by an 
$n$-simplex $\sigma_n$ in $K$, and it is its barycenter 
$\hat{\sigma}_n$ in $K'$. Since the vertices in $K'$ have the discrete 
topology, $\overline{e}^0_{\sigma_{n}}= e^0_{\sigma_{n}}$.
\item \label{ce1} If $d_1$ is a dual $1$-block, then there exists an 
$(n-1)$-simplex $\sigma_{n-1}$ in $K$ such that $d_1=D(\sigma_{n-1})$. 
Now, $\sigma_{n-1}$ is a face of exactly two $n$ simplices 
$\sigma^0_n$ and $\sigma_n^1$, and we have 
$$
\overline{d}_1=\overline{D}(\sigma_{n-1})= [\hat{\sigma}_n^0, 
\hat{\sigma}_{n-1} ] \cup  
 [\hat{\sigma}_n^1, \hat{\sigma}_{n-1}] \, .  
$$
(Notice that $\sigma_{n-1}$ could appear with any sign in
the boundaries $\partial \sigma_n^0$ and $\partial \sigma_n^1$.)
Thus, any $1$-cell $e^1_{\sigma_{n-1}}$ in $(K')^*$ is determined by
an $(n-1)$-simplex $\sigma_{n-1}$ in $K$, and it is given 
by the interior of a $1$-simplex of the 
form $[\hat{\sigma}_n, \hat{\sigma}_{n-1}]$ in $K'$,  
where $\sigma_n$ is an $n$-simplex in $K$ such that
$\sigma_{n} \succ \sigma_{n-1}$.
\item \label{ce2} If $d_2$ is a dual $2$-block, then there exists an
$(n-2)$-simplex in $K$ such that $d_2=D(\sigma_{n-2})$. Given any
$n$-simplex $\sigma_n$ that has $\sigma_{n-2}$ as a face, there exists
an $(n-1)$-simplex $\sigma_{n-1}$ such that
$$
\sigma_{n}\succ \sigma_{n-1} \succ \sigma_{n-2}\, , 
$$
and we have that 
$$
\overline{d}_2 =\overline{D}(\sigma_{n-2})=\cup_{\sigma_n, \sigma_{n-1}: \, 
\sigma_n\succ \sigma_{n-1}\succ 
\sigma_{n-2}} [\hat{\sigma}_n,\hat{\sigma}_{n-1}, \hat{\sigma}_{n-2}] 
$$
is the union of all such 2-simplices.
Thus, any $2$-cell $e^2_{\sigma_{n-2}}$ in $(K')^*$ is determined by an
$(n-2)$-simplex $\sigma_{n-2}$ in $K$, and it is given by the interior of 
a $2$-simplex of 
the form $[\hat{\sigma}_n, \hat{\sigma}_{n-1},\hat{\sigma}_{n-2}]$ in $K'$, 
where $\sigma_n$ and $\sigma_{n-1}$ are simplices in $K$ such that  
$\sigma_{n} \succ \sigma_{n-1} \succ \sigma_{n-2}$.
\item \label{ce3} If $d_3$ is a dual $3$-block, then there exists an 
$(n-3)$-simplex in $K$ such that $d_3=D(\sigma_{n-3})$. 
Given any $n$-simplex $\sigma_n$ that has $\sigma_{n-3}$ as a face, there are
simplices $\sigma_{n-1}$ and $\sigma_{n-2}$ of dimensions $n-1$ and $n-2$,
respectively, such that  
$$
\sigma_{n}\succ \sigma_{n-1} \succ \sigma_{n-2}\succ \sigma_{n-3}\, , 
$$
and we have that 
$$
\overline{d}_3=\overline{D}(\sigma_{n-3})=
\cup_{\sigma_n,\sigma_{n-1},\sigma_{n-2}:\, 
\sigma_{n}\succ \sigma_{n-1} \succ 
\sigma_{n-2}\succ \sigma_{n-3}}\, [\hat{\sigma}_n,\hat{\sigma}_{n-1},
\hat{\sigma}_{n-2},\hat{\sigma}_{n-3}] \, ,
$$
the union of all such $3$ simplices. 
Thus, any $3$-cell $e^3_{\sigma_{n-3}}$ in $(K')^*$ is determined by an 
$(n-3)$-simplex $\sigma_{n-3}$ in $K$, and it is given by the interior of
a $3$-simplex of the form $[\hat{\sigma}_n, \hat{\sigma}_{n-1},
\hat{\sigma}_{n-2},\hat{\sigma}_{n-3}]$ in $K'$, where $\sigma_n$, 
$\sigma_{n-1}$, and $\sigma_{n-2}$ are 
simplices in $K$ such that  
$\sigma_{n} \succ \sigma_{n-1} \succ \sigma_{n-2} \succ \sigma_{n-3}$.
\end{enumerate}

\subsection{The Stiefel-Whitney class $w_0(M)$}
Since every manifold is a $CW$ complex, and every $CW$ complex has the homotopy
type of a simplicial complex, it is natural to proceed by using 
cellular cohomology.

The Stiefel-Whitney class $w_0(M)$ is defined axiomatically
as the unit element $1 \in H^0(M; \mb{Z}/2)$ \cite[Axiom 1, p. 37]{mist}. 
Now, by \ref{ce0} above, a dual $0$-cell is a simplex of the form 
$e^0_{\sigma_n}=\hat{\sigma}_n$ in $K'$, $\sigma_n$ an $n$-simplex in $K$.
Thus, the $0$-cochain 
\begin{equation}
e^0_{\sigma_n}=\hat{\sigma}_n \mapsto 
w_0( \hat{\sigma}_n):= 1 \, {\rm mod}\; 2
\end{equation}
represents $w_0(M)$.

\section{Smoothly triangulated Riemannian manifolds: The classes 
$w_1(M)$, $w_2(M)$, and $w_3(M)$}
We endow the triangulated manifold manifold $M$ with an arbitrary
 Riemannian metric $g$. The Riemann curvature tensor of $g$ is  
$R^g(X,Y)Z=(\nabla^g_X\nabla^g_Y- \nabla^g_Y\nabla^g_X-
\nabla^g_{[X,Y]})Z$, where $\nabla^g$ the Levi-Civita connection of the metric.
It is usually expressed as the $(0,4)$ tensor
$g(R^g(X,Y)Z,W)$. The Ricci tensor $r_g(X,Y)$ of $g$ is 
the trace of the map $L\rightarrow R^g(L,X)Y$.
The scalar curvature $s_g$ is the metric trace of $r_g$. If
$(M,J,g)$ is Hermitian, the curvature $2$-form $\Omega^g(X,Y)$ of
the connection $\nabla^g$ is the trace of the map
$L\rightarrow R^g(L,JX)Y$.

Given any dual $i$-cell $e^{i}_{\sigma_{n-i}}$ in $(K')^*$ with closure
 $\overline{e}^i_{\sigma_{n-i}}=[\hat{\sigma}_n , \ldots, \hat{\sigma}_{n-i}]$,
$\sigma_n \succ \cdots \succ \sigma_{n-i}$, we fix its local orientation, and 
use the metric $g$ to regularize the corners and boundary
faces of the block $\overline{D}(\sigma_{n-i})$, as follows:

We let $\nu(\stackrel{\circ}{\sigma}_{n-i})$ denote the normal bundle of
$\stackrel{\circ}{\sigma}_{n-i} \subset \sigma_n$. We have
the local Whitney sum decomposition $T\sigma_n\mid_{\stackrel{\circ}{
\sigma}_{n-i}} = T\stackrel{\circ}{\sigma}_{n-i} \oplus
\nu (\stackrel{\circ}{\sigma}_{n-i})$. Since simplices are contractible, their 
tangent bundles are trivial. We choose the orientations of the fibers of these
three bundles compatibly, so that the product orientation of the fiber summands
in the Whitney sum equals the orientation of the fiber sum.  
Using the metric, we can then choose a positively oriented orthonormal
frame $\{e_1, \ldots , e_{n-i}\}$ for $T\stackrel{\circ}{\sigma}_{n-i}$,
and then extend it inductively to a positively oriented orthonormal frame 
$\{e_1, \ldots, e_n\}$ of $T\sigma_n\mid_{\stackrel{\circ}{\sigma}_{n-i}}$ 
such that, for each $j=1, \ldots, i$, $\{ e_1 , \ldots, e_{n-i+j}\}$ is a
positively oriented orthonormal frame for the tangent bundle 
$T\sigma_{n-i+j}\mid_{\stackrel{\circ}{\sigma}_{n-i}}$
over
$\stackrel{\circ}{\sigma}_{n-i}$ of the $(n-i+j)$th simplex in between.
By this construction, $\{ e_{n-i+1}, \ldots, e_n\}$ is a positively oriented
 orthonormal frame for the fibers of the normal bundle
$\nu(\stackrel{\circ}\sigma_{n-i})$. We say that 
$\{e_1, \ldots, e_n\}$ fixes compatibly the local positive orientation of 
the cell $e^i_{\sigma_{n-i}}$. When the manifold $M$ is oriented, we always 
choose the positive orientation of the starting 
$\{ e_1, \ldots, e_{n-i}\}$ in the construction to agree with the
orientation of $T\stackrel{\circ}{\sigma}_{n-i}$ induced by that of $M$,  
thus ensuring that the local positive orientation of 
$T\sigma_n \mid_{\stackrel{\circ}{\sigma}_{n-i}}$ as defined 
by the frame $\{e_1, \ldots, e_n\}$, and as defined by the 
orientation of $M$, then agree with each other.

We consider a tubular neighborhood of $\stackrel{\circ}{\sigma}_{n-i}$
in $M$. For some $\varepsilon > 0$, this neighborhood is obtained by applying 
the exponential map to vectors of norm less than $\varepsilon$ lying in 
$\nu_p(\stackrel{\circ}{\sigma}_{n-i})$, $p\in \stackrel{\circ}{\sigma}_{n-i}$ 
the fiber base. We choose $\varepsilon$ sufficiently small so that this action 
of the exponential map can be extended by continuity, with continuous inverse, 
to vectors of norm less or equal than $\varepsilon$. We denote by 
$E^\varepsilon(\stackrel{\circ}{\sigma}_{n-i})$ 
the resulting tubular neighborhood. It is the total space 
of a fiber bundle over $\stackrel{\circ}{\sigma}_{n-i}$  
whose fibers are geodesic open $i$-disks of radius $\varepsilon$ centered
at the base points. We denote the fiberwise closure  
of $E^\varepsilon(\stackrel{\circ}{\sigma}_{n-i})$ by
$\overline{E}^\varepsilon(\stackrel{\circ}{\sigma}_{n-i})$, and the fiberwise
boundary of this by 
$\partial \overline{E}^\varepsilon(\stackrel{\circ}{\sigma}_{n-i})$. 
They are the total spaces of fiber bundles over
$\stackrel{\circ}{\sigma}_{n-i}$ by closed geodesic $i$-disks,  
and $(i-1)$-spheres, centered at the base point, respectively.  

The fiber of $E^\varepsilon(\stackrel{\circ}{\sigma}_{n-i})$ through 
the barycenter
$\hat{\sigma}_{n-i}$, $D^\varepsilon_{\hat{\sigma}_{n-i}}$, is an open
$i$-disk with center at $\hat{\sigma}_{n-i}$, whose
closure $\overline{D}^\varepsilon_{\hat{\sigma}_{n-i}}$ is a manifold with 
boundary, the fiber
of $\overline{E}^\varepsilon(\stackrel{\circ}{\sigma}_{n-i})$ over
$\hat{\sigma}_{n-i}$. We have that 
$\partial \overline{D}^\varepsilon_{\hat{\sigma}_{n-i}}$ is the
fiber of $\partial \overline{E}^\varepsilon(\stackrel{\circ}{\sigma}_{n-i})$
over the said barycenter base point. Since 
$H_i(K,K\setminus \hat{\sigma}_{n-i}) \cong \mb{Z}$ for $i=n$, and vanishes 
otherwise, the block pair $(\overline{D}(\sigma_{n-i}),\dot{D}(\sigma_{n-i}))$ 
has the homology of an $i$-cell, modulo its boundary, and by construction, it 
is homotopically equivalent to the pair
$(\overline{D}^\varepsilon_{\hat{\sigma}_{n-i}}, \partial
\overline{D}^\varepsilon_{\hat{\sigma}_{n-i}})$. 
We have just smoothed out the corners and edges of the former. 
We provide this fiber $D^\varepsilon_{\hat{\sigma}_{n-i}}$, which is totally
geodesic in $M$, with the metric induced by $g$ on it that, for simplicity,
we denote by $g_{\hat{\sigma}_{n-i}}$. We call the pair 
$(D^\varepsilon_{\hat{\sigma}_{n-i}}, g_{\hat{\sigma}_{n-i}})$ the 
smooth $\varepsilon$-Riemannian $\sigma_{n-i}$ block. 
      
Under the inclusion map $i: \; \stackrel{\circ}{\sigma}_{n-i}\hookrightarrow 
M$, the pull-back bundle $i^* TM$ decomposes as
the Whitney sum $i^*TM = T\stackrel{\circ}{\sigma}_{n-1}\oplus
\nu(\stackrel{\circ}{\sigma}_{n-i})$. Since 
parallel transport along rays emanating from points of the zero section 
preserve inner products, the elements of the frame 
$\{e_1, \ldots, e_n\}$ of $T \sigma_n \mid_{\stackrel{\circ}{\sigma}_{n-i}}$ 
above, which fixes  compatibly the local positive orientation of the cell 
$e^i_{\sigma_{n-i}}$, can all be pushed forward under the exponential map 
from points on the said zero section to produce   
an orthonormal frame $\{ v_1, \ldots , v_n\}$ for $TM$ that is defined on
the tubular neighborhood $E^\varepsilon(\stackrel{\circ}{\sigma}_{n-i})$, and 
such that over $D^\varepsilon_{\hat{\sigma}_{n-i}}$, 
$\{ v_{n-i+1},\ldots, v_n\}$ is an orthonormal frame of
$TD^\varepsilon_{\hat{\sigma}_{n-i}}$. 
We say that this positively oriented orthonormal frame 
$\{ v_1, \ldots , v_n\}$ is compatible with the local positive orientation
of the smooth $\varepsilon$-Riemannian $\sigma_{n-i}$ block 
$(D^\varepsilon_{\hat{\sigma}_{n-i}}, g_{\hat{\sigma}_{n-i}})$. When 
$M$ is oriented, by construction, the positive orientation of 
the frame $\{ v_1, \ldots , v_n\}$, and that of 
$M$ agree with each another.

\subsection{The Stiefel-Whitney class $w_1(M)$} \label{w1}
By \S \ref{smt}\ref{ce1} above, a dual $1$-cell $e^1_{\sigma_{n-1}}$,  
determined by an $(n-1)$-simplex $\sigma_{n-1}$ in $K$, has closure given by
a simplex of the form 
$$
\overline{e}^1_{\sigma_{n-1}}=[\hat{\sigma}_n , \hat{\sigma}_{n-1}]
$$  
in $K'$, where $\sigma_n$ is an $n$-simplex in $K$ such that 
$\sigma_n \succ \sigma_{n-1}$. The dual $1$-block $D(\sigma_{n-1})$ 
that $\sigma_{n-1}$ determines is such that
$$
\overline{D}(\sigma_{n-1})=[\hat{\sigma}^0_n , \hat{\sigma}_{n-1}] \cup 
[\hat{\sigma}^1_n , \hat{\sigma}_{n-1}] \, ,
$$
where $\sigma_n^0$ and $\sigma_n^1$ are the two simplices in $K$ that
have $\sigma_{n-1}$ as a face. We choose any smooth $\varepsilon$-Riemannian 
$\sigma_{n-1}$ block $(D^{\varepsilon}_{\hat{\sigma}_{n-1}},
g_{\hat{\sigma}_{n-1}})$, and an  
orthonormal frame $\{v_1, \ldots , v_{n}\}$ for $TM\mid_{\stackrel{\circ}
{\sigma}_{n-1}}$, defined in a tubular neighborhood 
$E^\varepsilon(\stackrel{\circ}{\sigma}_{n-1})$ of the open simplex in $M$,  
whose orientation is compatible with that of the local positive orientation
of the block.
Thus, $\{v_1, \ldots, v_{n-1}\}$ are the push-forward under the exponential
map in the normal directions of the elements of a positively oriented 
orthonormal frame for 
$T\stackrel{\circ}{\sigma}_{n-1}$,
and $\{ v_n\}$ is a positively oriented frame for the fibers of the 
neighborhood of the zero section $\stackrel{\circ}{\sigma}_{n-1}$ of
$\nu(\stackrel{\circ}{\sigma}_{n-1})$ that lie in
$E^{\varepsilon}(\stackrel{\circ}{\sigma}_{n-1})$. By construction, on 
these points of the zero section of the bundle, the metric components are 
such that $g(v_i,v_j)=\delta_{i\, j}$, 
hence the function $\stackrel{\circ}{\sigma}_{n-1} \cap 
E^{\varepsilon}(\stackrel{\circ}{\sigma}_{n-1}) \ni p 
\rightarrow \det{g(v_i,v_j)}(p)\equiv 1$. In particular, 
\begin{equation} \label{sw1g}
e^1_{\sigma_{n-1}} \mapsto {\bf w}^g_1(e^1_{\sigma_{n-1}}) 
:= [\det{g(v_i, v_j)}(\hat{\sigma}_{n-1})] \;  
{\rm mod}\;   2 \, , 
\end{equation}
is a well defined $\mb{Z}/2$ $1$-cochain.       

\begin{theorem}
The cochain ${\bf w}_1^g$ in {\rm (\ref{sw1g})} represents the first
Stiefel-Whitney class $w_1(M)$. 
\end{theorem}

{\it Proof}. We have that $\det{g(v_i,v_j)}(p)=1$ for any
$p\in \stackrel{\circ}{\sigma}_{n-1} \cap E^{\varepsilon}(\stackrel{\circ}
{\sigma}_{n-1})
\subset \sigma_n\cap E^{\varepsilon}(\stackrel{\circ}{\sigma}_{n-1})$.
Therefore, 
$$
{\bf w}^g_1(e^1_{{\sigma}_{n-1}})= 1 \; {\rm mod}\;   2 \, , 
$$
and ${\bf w}^g_1$ is a $\mb{Z}/2$-cochain that assigns $1$ to each dual
$1$-cell in $(K')^*$. By Whitney's criteria \cite{whit} (see proof in 
\cite{chee}), this characterizes $w_1(M)$.
\qed

\subsection{The Stiefel-Whitney class $w_2(M)$} \label{w2}
By \S \ref{smt}\ref{ce2} above, a dual $2$-cell $e^2_{\sigma_{n-2}}$, 
determined by an $(n-2)$-simplex $\sigma_{n-2}$ in $K$, has closure that is
given by a simplex of the form 
$$
\overline{e}^2_{\sigma_{n-2}}=[\hat{\sigma}_n, \hat{\sigma}_{n-1},
\hat{\sigma}_{n-2}]
$$ 
in $K'$, where $\sigma_n$, and $\sigma_{n-1}$ are simplices in $K$ such
that $\sigma_{n} \succ \sigma_{n-1} \succ \sigma_{n-2}$.

We choose any smooth $\varepsilon$-Riemannian $\sigma_{n-2}$ block 
$(D^{\varepsilon}_{\hat{\sigma}_{n-2}}, g_{\hat{\sigma}_{n-2}})$, and 
an orthonormal frame $\{ v_1, \ldots, v_n\}$ for $TM\mid_{
\stackrel{\circ}{\sigma}_{n-2}}$, defined in a tubular neighborhood
$E^\varepsilon(\stackrel{\circ}{\sigma}_{n-2})$ of the open simplex in $M$,
and with orientation compatible with that of the local positive orientation
of the block. Thus, $D^\varepsilon_{\hat{\sigma}_{n-2}}$ is an open
$2$-disk, with center at $\hat{\sigma}_{n-2}$, whose
closure is a manifold with boundary, and carries the metric
$g_{\hat{\sigma}_{n-2}}$ induced by the ambient space metric $g$ on
it, which can be extended by continuity to a smooth metric on the closure.  
Since  $H_i(K,K\setminus \hat{\sigma}_{n-2})\cong \mb{Z}$ for
$i=n$, and vanishes otherwise,
the block pair $(\overline{D}(\sigma_{n-2}), \dot{D}(\sigma_{n-2}))$ has
the homology of a 2-cell modulo its boundary, and 
the smooth pair
$(\overline{D}^\varepsilon_{\hat{\sigma}_{n-2}}, \partial
\overline{D}^\varepsilon_{\hat{\sigma}_{n-2}})$ is homotopically equivalent to
it.
The orthonormal frame $\{ v_1, \ldots , v_n\}$ 
is such that over $D^\varepsilon_{\hat{\sigma}_{n-2}} \hookrightarrow 
E^\varepsilon(\stackrel{\circ}{\sigma}_{n-2})$,
$\{ v_{n-1}, v_n\}$ is an orthonormal frame of
$TD^\varepsilon_{\hat{\sigma}_{n-2}}$.
We define a complex structure $J_{\hat{\sigma}_{n-2}}$ on this $2$-disk 
by setting $J_{\hat{\sigma}_{n-2}}v_{n-1}:=v_n$. This complex structure 
induces the same local positive orientation on
the disk that it had already. 

By construction, $(D^\varepsilon_{\hat{\sigma}_{n-2}}, g_{\hat{\sigma}_{n-2}})$
is a totally geodesic submanifold of $(M,g)$. Hence, by Gauss' equation, 
the intrinsic and extrinsic Riemann curvature tensors are the same, and this
implies the relation  
\begin{equation} \label{ricci}
r_{g_{\hat{\sigma}_{n-2}}}(X,Y) =  
r_g(X,Y)-\sum_{i=1}^{n-2}g(R^g(v_i,X)Y,v_i) 
\end{equation}
between the intrinsic and extrinsic Ricci tensors (see \cite[\S 2]{gracie}). 
Since for dimensional reasons, we have that
$$
r_{g_{\hat{\sigma}_{n-2}}}(X,Y) = \frac{s_{g_{\hat{\sigma}_{n-2}}}}{2}
g(X,Y) \, ,
$$ 
by Gauss' theorem for geodesic triangles \cite{gauss},
and a limiting procedure that approximates small arcs by the geodesic segment
between their end points, we conclude that 
$$
\lim_{\varepsilon \searrow 0}\frac{1}{\varepsilon^2}
\int_{\overline{D}^\varepsilon_{\hat{\sigma}_{n-2}}} 
\frac{s_{g_{\hat{\sigma}_{n-2}}}}{2} g_{\hat{\sigma}_{n-2}}(
J_{\hat{\sigma}_{n-2}}v_{n-1},v_n)d\mu_{g_{\hat{\sigma}_{n-2}}} =
\lim_{\varepsilon \searrow 0}\frac{1}{\varepsilon^2}
\int_{\overline{D}^\varepsilon_{\hat{\sigma}_{n-2}}} 
\frac{s_{g_{\hat{\sigma}_{n-2}}}}{2}
d\mu_{g_{\hat{\sigma}_{n-2}}} = 2\pi \, . 
$$
Therefore, 
\begin{equation} \label{sw2g}
\begin{array}{rcl}
e^2_{{\sigma}_{n-2}} \mapsto {\bf w}^g_2(e^2_{{\sigma}_{n-2}}) & := & 
{\displaystyle \lim_{\varepsilon \searrow 0}\frac{1}{\varepsilon^2}\left(
 \frac{1}{2\pi}\int_{\overline{D}^\varepsilon_{
\hat{\sigma}_{n-2}}} \left( r_g(J_{\hat{\sigma}_{n-2}}v_{n-1},v_n)- 
 \right. \right.} \vspace{1mm} \\ & & \left. \left. \sum_{i=1}^{n-2}g(R^g(v_i,
J_{\hat{\sigma}_{n-2}}v_{n-1},v_n,v_i) \right)
d\mu_{g_{\hat{\sigma}_{n-2}}} \right)  \; {\rm mod}\;   2 
\end{array}
\end{equation}
is a well defined $\mb{Z}/2$ $2$-cochain.

\begin{theorem}
The cochain ${\bf w}_2^g$ in {\rm (\ref{sw2g})} represents the second 
Stiefel-Whitney class $w_2(M)$.
\end{theorem}

{\it Proof}. By construction, we have that
$$
\lim_{\varepsilon \searrow 0}
\frac{1}{2\pi \varepsilon^2}\int_{\overline{D}^\varepsilon_{\hat{\sigma}_{n-2}}}
(r_g(J_{\hat{\sigma}_{n-2}}v_{n-1},v_n)-\sum_{i=1}^{n-2}
g(R^g(v_i,J_{\hat{\sigma}_{n-2}}v_{n-1},v_n,v_i))
d\mu_{g_{\hat{\sigma}_{n-2}}} = 1 \, .
$$
Therefore, 
$$
e^2_{{\sigma}_{n-2}} 
\mapsto {\bf w}^g_2(e^2_{{\sigma}_{n-2}})
=1 \; {\rm mod}\;   2 \, ,  
$$
and ${\rm w}^g_2$ is a $\mb{Z}/2$-cochain that assigns value $1$ to each
$2$-cell in $(K')^*$. By Whitney's criteria \cite{whit} (see proof in 
\cite{chee}),   this characterizes $w_2(M)$.
\qed
\medskip

Notice that $(\overline{D}^\varepsilon_{\hat{\sigma}_{n-2}}, 
J_{\hat{\sigma}_{n-2}},g_{\hat{\sigma}_{n-2}})$ is K\"ahler, and that the 
curvature form of the Levi-Civita connection 
$\nabla^{g_{\hat{\sigma}_{n-2}}}$ is   
$$
\Omega^{g_{\hat{\sigma}_{n-2}}}= 
r_{g_{\hat{\sigma}_{n-2}}}(J_{\hat{\sigma}_{n-2}}v_{n-1},v_n)
d\mu_{g_{\hat{\sigma}_{n-2}}}
 = \frac{s_{g_{\hat{\sigma}_{n-2}}}}{2}d\mu_{g_{\hat{\sigma}_{n-2}}}\, .
$$
When $(M^{n=2m},J,g)$ is Hermitian, if 
$i: \overline{D}^\varepsilon_{\hat{\sigma}_{n-2}} \rightarrow M$ is the
inclusion map, by construction the tensors $i^*J$ and 
$J_{\hat{\sigma}_{n-2}}$, and  
$i^*d\mu_g$ and $d\mu_{g_{\hat{\sigma}_{n-2}}}$, coincide, respectively,
and we have that 
$$
i^*\Omega^g = (r_g(Jv_{n-1}, v_n)-\sum_{i=1}^{n-2}g(R^g(v_i,Jv_{n-1})v_n,v_i))
(i^*d\mu_{g}) \, .
$$
In this case, (\ref{ricci}) simply reads \cite[Eq. (3)]{gracie}  
\begin{equation}\label{riccigrac}
\Omega^{g_{\hat{\sigma}_{n-2}}}= 
i^*\Omega^g \, .
\end{equation}

\begin{corollary}
Suppose that $(M,J,g)$ is a Hermitian manifold. Then the cochain 
${\bf w}_2^g$ in {\rm (\ref{sw2g})} is the $\mb{Z}/2$ reduction of the 
first Chern class $c_1(M,J)$. 
\end{corollary}

{\it Proof}. As integral classes, we know that  
$$
c_1= \frac{1}{2\pi}[\Omega^g] \, .   
$$
Since (\ref{ricci}) is the fact that for totally geodesic $2$-submanifolds, 
identity (\ref{riccigrac}) between the curvature forms of the  
intrinsic and extrinsic Levi-Civita connections holds, the  
result follows by the ensuing definition of ${\bf w}_2^g$ derived from it.
\qed
\medskip

Chern classes were introduced by Chern in 1946 \cite{cher}, just a year
after he made seminal contributions to the then understanding of the 
Gauss-Bonnet theorem, the reason why this is known 
nowadays as the Chern-Gauss-Bonnet theorem. Notice that when 
$(M,J,g)$ is a K\"ahler manifold, if 
$\rho_g (X,Y)=r_g(JX,Y)$ is the Ricci form of the metric $g$,  
$\Omega^g= \rho_g$, and we have that $c_1=\frac{1}{2\pi}[\rho_g]$.   
The fact that in this case $(1/2\pi)\rho_g$ represents $c_1$ is basically 
proved in the alluded to paper, 
an excellent testimony to the power of Chern's ideas since the study of both, 
complex and K\"ahler structures, were then in their infancies.
 
\subsection{The Stiefel-Whitney class $w_3(M)$} \label{w3}
By \S \ref{smt}\ref{ce3} above, a dual $3$-cell $e^3_{\sigma_{n-3}}$, 
determined by an $(n-3)$-simplex $\sigma_{n-3}$ in $K$, has as closure a 
simplex of the form 
$$
\overline{e}^3_{\sigma_{n-3}}=[ \hat{\sigma}_n, \hat{\sigma}_{n-1},\hat{\sigma}_{n-2},\sigma_{n-3}]
$$
in $K'$, where $\sigma_n$, $\sigma_{n-1}$, and $\sigma_{n-2}$ are
simplices in $K$ such that 
$\sigma_{n} \succ \sigma_{n-1} \succ \sigma_{n-2}\succ \sigma_{n-3}$.

We choose any smooth $\varepsilon$-Riemannian $\sigma_{n-3}$ block 
$(D^{\varepsilon}_{\hat{\sigma}_{n-3}}, g_{\hat{\sigma}_{n-3}})$, and 
an orthonormal frame $\{ v_1, \ldots, v_n\}$ for $TM\mid_{
\stackrel{\circ}{\sigma}_{n-3}}$, defined in a tubular neighborhood
$E^\varepsilon(\stackrel{\circ}{\sigma}_{n-3})$ of the open simplex in $M$,
and with orientation compatible with that of the local positive orientation
of the smooth block. Thus, $D^\varepsilon_{\hat{\sigma}_{n-3}}$, is an open
totally geodesic $3$-disk, with center at $\hat{\sigma}_{n-3}$, whose
closure is a manifold with boundary, and the orthonormal frame
is such that over $D^\varepsilon_{\hat{\sigma}_{n-3}}$, 
$\{ v_{n-2},v_{n-1},v_n\}$ is an orthonormal frame of
$TD^\varepsilon_{\hat{\sigma}_{n-3}}$.
We have that $H_i(K,K\setminus \hat{\sigma}_{n-3})\cong \mb{Z}$ for
$i=n$, and vanishes otherwise, so the block
$(\overline{D}(\sigma_{n-3}), \dot{D}(\sigma_{n-3}))$ has
the homology of a 3-cell modulo its boundary, and this pair and the
smooth pair $(\overline{D}^\varepsilon_{\hat{\sigma}_{n-3}}, \partial
\overline{D}^\varepsilon_{\hat{\sigma}_{n-3}})$ are homotopically equivalent.

By continuity, we extend the metric $g_{\hat{\sigma}_{n-3}}$ on 
the oriented geodesic $3$-disk $D^{\varepsilon}_{\hat{\sigma}_{n-3}}$ to a
metric on $\overline{D}^{\varepsilon}_{\hat{\sigma}_{n-3}}$. It
induces a volume form $d\mu_{g_{\hat{\sigma}_{n-3}}}$ on the oriented 
closed disk,
and a compatibly oriented area form $d\sigma_{g_{\hat{\sigma}_{n-3}}}$ on 
the geodesic $2$-sphere 
$\partial \overline{D}^\varepsilon_{\hat{\sigma}_{n-3}}$.
We may approximate the $3$-disk 
$\overline{D}^{\varepsilon}_{\hat{\sigma}_{n-3}}$
by a polyhedron centered at ${\hat{\sigma}_{n-3}}$
whose boundary polygons have total area yielding an approximation to the
area of  $\partial \overline{D}^{\varepsilon}_{\hat{\sigma}_{n-3}}$, of
quadratic order in $\varepsilon$. Then,
$$
\lim_{\varepsilon \searrow 0}\frac{1}{\varepsilon ^2} 
\int_{\partial \overline{D}^\varepsilon_{\hat{\sigma}_{n-3}}} 
d\sigma_{g_{\hat{\sigma}_{n-3}}}  =4\pi \, , 
$$
(the total solid angle subtended by 
$\partial \overline{D}^\varepsilon_{\hat{\sigma}_{n-3}}$), and we have that
\begin{equation} \label{sw3g}
e^3_{\sigma_{n-3}} 
 \mapsto
{\bf w}^g_3(e^3_{\sigma_{n-3}}) := 
\lim_{\varepsilon \searrow 0}\frac{1}{4\pi \varepsilon ^2}       
\int_{\partial \overline{D}^\varepsilon_{\hat{\sigma}_{n-3}}} 
d\sigma_{g_{\hat{\sigma}_{n-3}}} \; {\rm mod} \; 2 
\end{equation}
is a well defined $\mb{Z}/2$ $3$-cochain.

\begin{theorem}
The cochain ${\bf w}^g_3$ in {\rm (\ref{sw3g})} represents the third 
Stiefel-Whitney class $w_3(M)$.  
\end{theorem}

{\it Proof}. By construction, 
 ${\bf w}^g_3(e^3_{\sigma_{n-3}}) =1$ mod $2$,
so ${\bf w}^g_3$ is $\mb{Z}/2$-cochain that assigns value $1$ to each
$3$-cell in $(K')^*$. By Whitney's criteria  
\cite{whit} (see proof in \cite{chee}), this characterizes $w_3(M)$.  
\qed

\section{Metric representatives of higher degree Stiefel-Whitney classes}
The argument in \S \ref{w2} for $w_2(M)$ generalizes to produce 
representatives of the even degree classes $w_{2k}(M)$, $k>1$. The argument
in \S \ref{w3} for $w_3(M)$ generalizes to produce representatives of
the odd degrees classes $w_{2k+1}(M)$, $k>1$. We prove these assertions 
here.  

For convenience, we denote by $\omega_{d}$ the volume of the unit 
$d$-sphere $\mb{S}^d$ in Euclidean space $\mb{R}^{d+1}$.  

\subsection{The Stiefel-Whitney class $w_{2k}(M)$}
A dual $2k$-cell $e^{2k}_{\sigma_{n-2k}}$, determined by an 
$(n-2k)$-simplex $\sigma_{n-2k}$ of
$K$, has closure a $2k$-simplex of the form
$$
\overline{e}^{2k}_{\sigma_{n-2k}}=[\hat{\sigma}_{n}, \ldots , 
\hat{\sigma}_{n-2k}]
$$
in $K'$, where $\sigma_{l+1}\succ \sigma_{l}$, $l=n-2k, \ldots, n-1$.

We choose any smooth $\varepsilon$-Riemannian $\sigma_{n-2k}$ block 
$(D^{\varepsilon}_{\hat{\sigma}_{n-2k}}, g_{\hat{\sigma}_{n-2k}})$, and 
an orthonormal frame $\{ v_1, \ldots, v_n\}$ for $TM\mid_{
\stackrel{\circ}{\sigma}_{n-2k}}$, defined in a tubular neighborhood
of the open simplex in $M$ containing 
$D^{\varepsilon}_{\hat{\sigma}_{n-2k}}$,
 with orientation compatible with that of the 
local positive orientation of the smooth block. 
Then, $D^\varepsilon_{\hat{\sigma}_{n-2k}}$ is an open oriented 
$2k$-disk, with center at $\hat{\sigma}_{n-2k}$, whose
closure is a manifold with boundary.
We have that $H_i(K,K\setminus \hat{\sigma}_{n-2k})\cong \mb{Z}$ for
$i=n$, and vanishes otherwise, so the block pair 
$(\overline{D}(\sigma_{n-2k}), \dot{D}(\sigma_{n-2k}))$ has
the homology of a $2k$-cell modulo its boundary,
and the smooth pair
$(\overline{D}^\varepsilon_{\hat{\sigma}_{n-2k}}, \partial
\overline{D}^\varepsilon_{\hat{\sigma}_{n-2k}})$ is homotopically equivalent 
to it. 
The orthonormal frame $\{ v_1, \ldots , v_n\}$ is such that, over
$D^\varepsilon_{\hat{\sigma}_{n-2k}}$,
$\{ v_{n-2k+l}\}_{l=1}^{2k}$ is a positively oriented orthonormal frame of
$TD^\varepsilon_{\hat{\sigma}_{n-2k}}$.
We define a complex structure $J_{\hat{\sigma}_{n-2k}}$ on this $2k$-disk
by setting $J_{\hat{\sigma}_{n-2}}v_{n-2k+(2j+1)}:=v_{n-2k+(2j+2)}$, $j=0, 
\ldots, k-1$. 
This complex structure induces the same local positive orientation on
the disk that it had already.

By construction, $(D^\varepsilon_{\hat{\sigma}_{n-2k}}, 
g_{\hat{\sigma}_{n-2k}})$
is a totally geodesic submanifold of $(M,g)$, and with the metric extended
continuously to the closed disk, 
$(\overline{D}^\varepsilon_{\hat{\sigma}_{n-2k}},J_{\hat{\sigma}_{n-2k}},
g_{\hat{\sigma}_{n-2k}})$ is Hermitian. We let 
$\Omega^{g_{\hat{\sigma}_{n-2k}}}$ be the
curvature $2$-form of the Levi-Civita connection 
$\nabla^{g_{\hat{\sigma}_{n-2k}}}$. Then, by the Chern-Gauss-Bonnet theorem
for polyhedral manifolds of Allendoerfer and Weil \cite{alwe},  
we conclude that
$$
\frac{2^k}{(2k)!\varepsilon^{2k}}
\int_{\overline{D}^\varepsilon_{\hat{\sigma}_{n-2k}}}
\Omega^{g_{\hat{\sigma}_{n-2k}}} \wedge \cdots \wedge
\Omega^{g_{\hat{\sigma}_{n-2k}}} \rightarrow \frac{(2\pi)^k}{(2k-1)(2k-3)\cdots
3 \cdot 1}=\frac{\omega_{2k}}{2} \, ,  
$$
and therefore, the expression 
\begin{equation} \label{sw2kg}
e^{2k}_{{\sigma}_{n-2k}} \mapsto {\bf w}^g_2(e^{2k}_{{\sigma}_{n-2k}}):=
\lim_{\varepsilon \searrow 0}\frac{2^{k+1}}{\omega_{2k}
(2k)!\varepsilon^{2k}}\int_{\overline{D}^\varepsilon_{\hat{\sigma}_{n-2k}}}
\Omega^{g_{\hat{\sigma}_{n-2k}}} \wedge \cdots \wedge
\Omega^{g_{\hat{\sigma}_{n-2k}}} 
 \; {\rm mod}\;   2
\end{equation}
is a well defined $\mb{Z}/2$ $2k$-cochain that assigns the value $1\in 
\mb{Z}/2$ to the
dual cell $e^{2k}_{{\sigma}_{n-2k}}$.

\begin{theorem}
The cochain ${\bf w}^g_{2k}$ in {\rm (\ref{sw2kg})} represents the $2k$th
Stiefel-Whitney class $w_{2k}(M)$.
\end{theorem}

If $(M^{n=2m},J,g)$ is Hermitian, and 
$i: \overline{D}^\varepsilon_{\hat{\sigma}_{n-2k}} \rightarrow M$ is the
inclusion map, by construction we have that the tensors $i^*J$ and
$J_{\hat{\sigma}_{n-2k}}$, and
$i^*d\mu_g$ and $d\mu_{g_{\hat{\sigma}_{n-2k}}}$, coincide, respectively,
and we have that 
$i^*\Omega^g =\Omega^{g_{\hat{\sigma}_{n-2}}}$. Thus, the identity
\begin{equation}\label{grac}
\Omega^{g_{\hat{\sigma}_{n-2k}}}\wedge \cdots \wedge 
\Omega^{g_{\hat{\sigma}_{n-2k}}} =
 (i^*\Omega^g)\wedge \cdots \wedge (i^*\Omega^g)
\end{equation}
of the $k$-fold products of the intrinsic and extrinsic curvature forms
holds. 

\begin{corollary}
Suppose that $(M,J,g)$ is Hermitian.
Then the cochain ${\bf w}_{2k}^g$ in
{\rm (\ref{sw2kg})} is the $\mb{Z}/2$ reduction of the $k$th Chern
class $c_k(M,J)$.
\end{corollary}

\subsection{The Stiefel-Whitney class $w_{2k+1}(M)$}
A dual $(2k+1)$-cell $e^{2k+1}_{\sigma_{n-2k}}$, determined by an 
$(n-2k-1)$-simplex $\sigma_{n-2k-1}$ of $K$, has closure a $(2k+1)$-simplex 
of the form
$$ 
\overline{e}^{2k+1}_{\sigma_{n-2k-1}}=[\hat{\sigma}_{n},\hat{\sigma}_{n-1}, 
\ldots , \hat{\sigma}_{n-2k-1}]
$$ 
in $K'$, where any pair $\sigma_{l+1}, \sigma_l$ of consecutive 
intermediate simplices are such that $\sigma_{l+1}\succ \sigma_{l}$. 

We choose any smooth $\varepsilon$-Riemannian $\sigma_{n-2k-1}$ block
$(D^{\varepsilon}_{\hat{\sigma}_{n-2k-1}}, g_{\hat{\sigma}_{n-2k-1}})$.
Thus, $D^\varepsilon_{\hat{\sigma}_{n-2k-1}}$, is a totally geodesic
open $(2k+1)$-disk, with center at $\hat{\sigma}_{n-2k-1}$, and whose
closure is a manifold with boundary. 
We have that $H_i(K,K\setminus \hat{\sigma}_{n-2k-1})\cong \mb{Z}$ for
$i=n$, and vanishes otherwise, so the block
pair $(\overline{D}(\sigma_{n-2k-1}), 
\dot{D}(\sigma_{n-2k-1}))$ has the homology of a $(2k+1)$-cell modulo
its boundary, and the smooth pair
$(\overline{D}^\varepsilon_{\hat{\sigma}_{n-2k-1}}, \partial
\overline{D}^\varepsilon_{\hat{\sigma}_{n-2k-1}})$ is homotopically 
equivalent to it.

With the continuously extended metric $g_{\hat{\sigma}_{n-2k-1}}$, the 
pair $(\overline{D}^{\varepsilon}_{\hat{\sigma}_{n-2k-1}},
g_{\hat{\sigma}_{n-2k-1}})$ is a totally geodesic $(2k+1)$ submanifold 
with boundary of the ambient $(M,g)$. The metric
induces a volume form $d\mu_{g_{\hat{\sigma}_{n-2k-1}}}$ on it,
and a compatibly oriented hypersurface area form 
$d\sigma_{g_{\hat{\sigma}_{n-2k-1}}}$ on   
the $2k$-sphere $\partial \overline{D}^\varepsilon_{\hat{\sigma}_{n-2k-1}}$.
Since we have  
$$
\lim_{\varepsilon \searrow 0}\frac{1}{\varepsilon ^{2k}} 
\int_{\partial \overline{D}^\varepsilon_{\hat{\sigma}_{n-2k-1}}} 
d\sigma_{g_{\hat{\sigma}_{n-2k-1}}} =
\frac{2^{k+1}\pi^k}{(2k-1)(2k-3)\cdots 3 \cdot 1}=
\omega_{2k}\, , 
$$
(the total solid angle subtended by 
$\partial \overline{D}^\varepsilon_{\hat{\sigma}_{n-2k-1}}$), the 
expression 
\begin{equation} \label{sw3kg}
e^{2k+1}_{\sigma_{n-2k-1}} \mapsto 
{\bf w}^g_{2k+1}{(e^{2k+1}_{\sigma_{n-2k+1}})} := 
\lim_{\varepsilon \searrow 0}\frac{1}{\omega_{2k} \varepsilon ^{2k}}       
\int_{\partial \overline{D}^\varepsilon_{\hat{\sigma}_{n-2k-1}}} 
d\sigma_{g_{\hat{\sigma}_{n-2k-1}}} \; {\rm mod} \; 2 
\end{equation}
is a well defined $\mb{Z}/2$ $(2k+1)$-cochain assigning the value $1\in 
\mb{Z}/2$ to the dual cell $e^{2k+1}_{\sigma_{n-2k-1}}$. 

\begin{theorem} 
The cochain ${\bf w}^g_{2k+1}$ in {\rm (\ref{sw3kg})} represents the 
$(2k+1)$th Stiefel-Whitney class $w_{2k+1}(M)$.
\end{theorem}

All the odd dimensional dual cells come in pairs: Given 
the cell $e^{2k+1}_{\sigma_{n-2k-1}}$ associated to the chain of 
simplices $\sigma_n \succ \sigma_{n-1} \succ \cdots \succ \sigma_{n-2k-1}$,
there is exactly one other $(2k+1)$-cell $\tilde{e}^{2k+1}_{\sigma_{n-2k-1}}$ 
determined by $\sigma_{n-2k-1}$ that is associated to a chain of simplices of 
the form
$\tilde{\sigma}_n \succ \sigma_{n-1} \succ \cdots \succ \sigma_{n-2k-1}$,
where all but the first of the simplices in the latter chain are the same as 
those in the former one. For the $(n-1)$-simplex $\sigma_{n-1}$ in the 
first chain is a face of exactly two simplices of top dimension, 
$\sigma_n$, and a second one, which we call $\tilde{\sigma}_n$, and use it 
to form the second chain. Notice that if
$(M,K)$ is an oriented homology $n$-manifold,   
and we orient all the $n$-simplices in $K$ so that 
$\gamma = \sum \sigma_n $ is the cycle $1$ in $H_n(M;\mb{Z})$, 
and orient the other simplices of $K$ arbitrarily, then   
the orientations of $\sigma_n$ and $\tilde{\sigma}_n$ are such that 
the boundary $(n-1)$-chain $\partial \sigma_n + \partial \tilde{\sigma}_n$ has 
coefficient zero on $\sigma_{n-1}$. Thus, if $\partial \sigma_n$ has 
coefficient $1$ on $\sigma_{n-1}$, then $\partial \tilde{\sigma}_n$ has
coefficient $-1$ on it, and vice versa. (Naturally, when $M$ is an oriented
manifold, the orientation of $M$ as a homological $n$-manifold that we use, 
and the orientation of $M$ as a manifold, agree with each other.) 

Let us assume that $(M^{n=2m},J,g)$ is Hermitian.
If $k=0$, since ${\bf w}_{1}^g$ maps any dual $1$-cell to $1$, and the local 
positive orientations of simplices are all compatible with the orientation
of $M$, 
we conclude that the local compatible positive orientations of 
$e^{1}_{\sigma_{n-1}}$ and $\tilde{e}^{1}_{\sigma_{n-1}}$ must be, in turn,  
compatible with each other, and since this conclusion is independent of
$\sigma_{n-1}$, compatible with the local positive 
orientation of any $(n-1)$-dual cell. Therefore, 
${\bf w}_{1}^g=1$ in cohomology.
By an induction on $k$, working on the $(n-2k)$-skeleton at the time (skeleton
on which $J$ induces a natural orientation compatible with that of $M$ as a 
whole), 
we handle all the possible choices of
intermediate simplices larger than $\sigma_{n-2k+1}$ in the chain 
$\sigma_n \succ \sigma_{n-1} \succ \cdots \succ \sigma_{n-2k-1}$
for a given $\sigma_{n-2k-1}$, and conclude that 
${\bf w}^g_{2k+1}=1$, so ${\bf w}^g_{2k+1}$ is then trivial in cohomology.

\begin{theorem}
Suppose that $(M,J,g)$ is Hermitian. Then the cochain ${\bf w}^g_{2k+1}$ in
{\rm (\ref{sw3kg})} is trivial in cohomology, and so any odd degree 
Stiefel-Whitney class $w_{2k+1}(M)$ is trivial.
\end{theorem}

Notice that the starting step in the induction argument above is really 
dependent
on the orientation of $M$ only, so if we apply it to the cochain 
${\bf w}^g_1$ in (\ref{sw1g}) for oriented manifolds of any dimension, we
would conclude then that ${\bf w}^g_1$ is trivial in cohomology, and therefore, 
so would be $w_{1}(M)$, as expected.


\begin{thebibliography}{CA}
\bibitem{alwe}
C.B. Allendoerfer \& A. Weil, {\it The Gauss-Bonnet Theorem for 
Riemannian Polyhedra}, Trans. Amer. Math. Soc., 53 (1943), pp. 101-129. 
\bibitem{chee}
J. Cheeger, {\it A combinatorial formula for Stiefel-Whitney classes},
Topology of Manifolds, J.C. Cantrell \& C.H. Edwards eds., Markham Publ.
Co. (1970), pp. 470-471.
\bibitem{cher}
S. Chern, {\it Characteristic classes of Hermitian manifolds}, Ann. of Math., 
47 (1946), pp. 85-121.
\bibitem{gauss}
C.F. Gauss, {\it Disquisitiones generales circa superficies
curvas}, Commentationes Gottingensis, 1827.
\bibitem{hato}
S. Halperin \& D. Toledo, {\it Stiefel-Whitney homology classes},
Ann. of Math., 96 (1972), pp. 511-525.
\bibitem{hazw}
R. Harvey \& J. Zweck, {\it Stiefel-Whitney currents}, J. Geom. Anal.
8 (1998), pp. 809-844. 
\bibitem{mist}
J.W. Milnor \& J.D. Stasheff, {\it Characteristic Classes}. Annals of
Mathematics Studies 76, Pricenton University Press, 1974.
\bibitem{mu}
J.R. Munkres, {\it Concordance is equivalent to smoothability}, Topology 
5 (1966) pp. 371-389.
\bibitem{mubo}
J.R. Munkres, {\it Elements of algebraic topology}. Addison-Wesley
Publishing Company, Menlo Park, CA, 1984. ix+454 pp.
\bibitem{gracie}
S.R. Simanca, {\it Isometric Embeddings I: General Theory}. Riv. Mat.
Univ. Parma, 8 (2017), pp. 307-343.
\bibitem{stie}
E. Stiefeld, {\it Richtungsfelderund Fernparallelismus in 
Mannigfaltigkeiten}, Comm. Math. Helv.  8 (1936), pp. 
\bibitem{sull}
D. Sullivan, {\it Combinatorial invariants of analytic spaces}. 1971 
Proceedings of Liverpool Singularities-Symposium, I (1969/70),
pp. 165–168, Springer, Berlin. 
\bibitem{whit}
H. Whitney, {\it On the theory of sphere-bundles}, Proc. Nat. Acad. 
Sci. U.S.A. 26 (1940), pp. 148–153.
\end{thebibliography}
\end{document}